\documentclass[a4paper,12pt,onecolumn]{article}

\usepackage{hyperref}
\usepackage[vmargin=2cm,hmargin=2cm,headheight=14.5pt,top=2cm,headsep=.5cm]{geometry}

\usepackage{bm}
\usepackage[utf8]{inputenc}

\usepackage{empheq}
\usepackage{stackrel}
\usepackage{cases}
\usepackage{mathtools}
\usepackage{amsthm,amsmath,amscd}
\usepackage{makeidx}
\usepackage[charter]{mathdesign}
\usepackage{perpage}
\usepackage{bm}
\usepackage{tikz-cd}
\tikzcdset{every label/.append style = {font = \small}}
\usepackage{caption}
\usepackage{graphicx}
\DeclareMathSizes{12}{12}{8}{6}

\usepackage{cite}
\usepackage{url}
\usepackage[ddmmyy]{datetime}
\usepackage{enumitem}
\usepackage{framed}
\usepackage{arydshln}
\usepackage[symbol]{footmisc}
\usepackage{marvosym}
\usepackage[scaled=0.85]{beramono}
\usepackage{perpage}
\MakePerPage[1]{footnote}

\makeatletter
\renewcommand*{\@makefnmark}{\hbox{\@textsuperscript{%
			\normalfont\@thefnmark}}}
\makeatother
\makeatletter
\def\@fnsymbol#1{\ensuremath{\ifcase#1\or \text{\Mercury}
		\or \text{\Venus} \or \text{\Earth} \or
		\text{\Jupiter} \or \text{\Saturn} \or \text{\Neptune} \or \text{\Uranus} \or
		\text{\Pluto}
		\or \text{\Moon} \or \text{\Sun}
		\else\@ctrerr\fi}}%
\makeatother

\newtheoremstyle{ptheorem}{1em}{0em}{\itshape}{}{\bfseries}{.}{.5em}{\thmname{#1}\thmnumber{
		#2}\thmnote{ (\hspace{-.01pt}{#3})}}

\theoremstyle{ptheorem}

\newtheorem{thm}{Theorem}[section]

\newtheorem{lem}[thm]{Lemma}
\newtheorem{cor}[thm]{Corollary}

\newtheoremstyle{hdef}{1em}{0em}{}{}{\bfseries}{.}{.5em}{\thmname{#1}\thmnumber{
		#2}\thmnote{ (\hspace{-.01pt}{#3})}}
\theoremstyle{hdef}

\newtheorem{dfn}[thm]{Definition}
\newtheorem{rem}[thm]{Remark}

\newtheorem{exa}[thm]{Example}

\numberwithin{equation}{section}
\numberwithin{figure}{section}

\let\originalleft\left
\let\originalright\right
\renewcommand{\left}{\mathopen{}\mathclose\bgroup\originalleft}
\renewcommand{\right}{\aftergroup\egroup\originalright}

\DeclareMathOperator{\dif}{d}

\DeclareMathOperator{\conv}{conv}

\newcommand{\cB}{{\mathcal B}}
\newcommand{\cC}{{\mathcal C}}

\newcommand{\cL}{{\mathcal L}}
\newcommand{\cM}{{\mathcal M}}
\newcommand{\cN}{{\mathcal N}}

\newcommand{\bN}{{\mathbb N}}

\newcommand{\bR}{{\mathbb R}}

\renewcommand{\l}{\lambda}
\newcommand{\e}{\varepsilon}
\renewcommand{\L}{\Lambda}

\renewcommand{\phi}{\varphi}
\renewcommand{\le}{\leqslant}
\renewcommand{\ge}{\geqslant}

\newcommand{\ol}{\overline}

\newcommand{\n}{{n\in\bN}}

\newcommand{\Ra}{\Rightarrow}

\renewcommand{\d}{\delta}

\renewcommand{\(}{\left(}
\renewcommand{\)}{\right)}

\newcommand{\til}{\widetilde}
\newcommand{\pd}{\partial}

\newcommand{\bs}{\backslash}

\newcommand{\olb}[1]{%
	\vbox{\offinterlineskip\ialign{\hfil##\hfil\cr
			$\rotatebox[origin=c]{90}{$]$}$\cr\noalign{\kern-.45ex}{$#1$}\cr}}}

\newcommand{\noop}[1]{}

\renewcommand{\ss}{\subset}

\allowdisplaybreaks

\begin{document}

\title{When is the composition of functions measurable?}

\date{}

\author{
F. Javier Fernández\\
	\normalsize e-mail: fjavier.fernandez@usc.es\\
	F. Adri\'an F. Tojo \\
	\normalsize e-mail: fernandoadrian.fernandez@usc.es\\
	\normalsize \emph{CITMAga, 15782, Santiago de Compostela, Spain},\\ \normalsize \emph{Departamento de Estatística, Análise Matemática e Optimización},\\ \normalsize \emph{Universidade de Santiago de Compostela, 15782, Facultade de Matemáticas, Santiago, Spain.}}

 \maketitle
 
\begin{abstract}
	In this article we explore under which conditions on the interior function the composition of functions is measurable. We also study the sharpness of the result by providing a counterexample for weaker hypotheses.
\end{abstract}


 \section{Introduction}

It is a well known fact among young analysts that \emph{the composition of measurable functions is not necessarily measurable}, although this fact, to be true, has to be stated precisely, as it depends on what we mean when we say that a function is measurable. In this case, we are referring to the fact that the inverse image of an open set is a Lebesgue-measurable set of the real line. It is therefore convenient to  start with the definition of measurable function in a broader sense.
\begin{dfn}
	Let $(X,\cM)$, $(X,\cN)$ be measurable spaces. We say $f:(X,\cM)\to(X,\cN)$ is measurable if $f^{-1}(A)\in\cM$ for every $A\in \cN$.
\end{dfn}
It is clear then, from this definition, that whether a function is measurable or not depends on the measurable spaces chosen. 
Let $\cL$ and $\cB$ be the Lebesgue and Borel $\sigma$-algebras on $\bR$ respectively. In general, if $f,g:(\bR,\cL)\to(\bR,\cB)$ are measurable $g\circ f:(\bR,\cL)\to(\bR,\cB)$ does not have to be so, which is the case we first talked about. In order for $g\circ f:(\bR,\cL)\to(\bR,\cB)$ to be measurable, it would be enough if  $f: (\bR,\cL)\to(\bR,\cL)$ or $g: (\bR,\cB)\to(\bR,\cB)$ were measurable (this last condition holds when $g$ is continuous), but in general this is not the case.

Here we are interested in sufficient conditions on $f:(\bR,\cL)\to(\bR,\cB)$  that can guarantee that, if $g:(\bR,\cL)\to(\bR,\cB)$ is measurable, so is $g\circ f:(\bR,\cL)\to(\bR,\cB)$. We will consider then $\mu^*$ and $\mu$, the Lebesgue exterior measure and  measure respectively. For brevity, we will speak of $\cL/\cB$ and $\cL/\cL$ measurable functions.

Rather counterintuitively, great regularity or monotonicity of $f$ does not guarantee that the composition $g\circ f$ will be $\cL/\cB$ measurable when $g$ is. We illustrate this point with the following example.

\begin{exa}\label{exace}
	We will first proceed to construct a strictly increasing $\cC^\infty$ function $f$  which is not $\cL/\cL$ measurable, that is, such that there exists $D\in\cL$ satisfying $f^{-1}(D)\not\in \cL$.

We start by recalling the construction of the Smith–Volterra–Cantor set (or fat Cantor set) \cite[p.~39]{Folland}. Let $\l\in(0,1)$ and define $r:=\frac{1-\l}{3-2 \l}$. Observe that $r\in(0,1/3)$. 
The \emph{Smith–Volterra–Cantor set $C$ of measure $\l$} is created  iteratively in the following way. Let $C_0 := [0,1]$ and, given $C_{n-1}$ define $C_n$  by deleting the central open interval of length $r^n$ from  each the connected components of $C_{n-1}$. It is clear by construction that $C_n\ss C_{n-1}$ and $C_n$ is compact for every $n\ge 1$. Thus, the set $C:= \bigcap_{n=0}^\infty C_n$, the Smith–Volterra–Cantor set, is well defined, non-empty and compact. Since, by definition, $C_n$ has $2^{n-1}$ connected components, \[\mu(C_n)=\mu(C_{n-1})-2^{n-1}r^n=1-\frac{r \left(1-2^n r^n\right)}{1-2 r}\]
for every $n\in\bN$, we conclude that \[\mu(C)=\lim_{n\to\infty}\mu(C_n)=1-\frac{r}{1-2 r}=\frac{1-3r}{1-2 r}=\l.\]

    Let $C\ss I:=[0,1]$ be a Smith–Volterra–Cantor set such that $\mu(C)\in(0,1)$. The positive measure of this set will be crucial for this example --cf. Theorem~\ref{dna1}. Consider the function
    \[\psi(x):=e^{-(1-x^2)^{-1}}\]
for $x\in(-1,1)$. $\psi\in\cC^\infty((-1,1),\bR)$. Since $C$ is closed, $I\bs C$ is open, so it has a countable number of connected components each of which is an open interval. Let us define $h(x)=0$ if $x\in C$ and \[h(x)=2^{-(b-a)^{-1}}\psi\(\frac{2x-b-a}{b-a}\)\]
 if $x\in(a,b)$ where $(a,b)$ is a connected component of $I\bs C$. Let $M_n:=\max \lvert\psi^{(n)}\rvert$. We will show now that $h\in\cC^\infty(I,\bR)$. It is clear that $h$ is $\cC^\infty$ in  $I\bs C$. If $x\in C$, let us check that $\lim_{y\to x^-}f(x)=0$ (in case the limit can be taken). If $x=b$ for some $(a,b)\ss I\bs C$, this is obvious. Otherwise, there is a sequence of points in $C$ converging to $x$ from the left. Thus, given $\e\in\bR^+$, there exists $y\in C$, $x-\e<y<x$, so any $z\in (y,x)\bs C$ belongs to an interval $(a,b)$ with $b-a<\e$ and, therefore, \[h(z)=2^{-(b-a)^{-1}}\psi((2x-b-a)/(b-a))<2^{-(b-a)^{-1}}<b-a<\e.\] Hence,  $\lim_{y\to x^-}f(x)=0$. Repeating the argument for limits from the right and observing that $h|_C=0$, we conclude that $h$ is continuous. Assuming $h$ is $n-1$ times differentiable and taking into account that
\[\lvert h^{(n)}(x)\rvert=2^{-(b-a)^{-1}}2^n(b-a)^{-n}\left\lvert\psi^{(n)}\(\frac{2x-b-a}{b-a}\)\right\rvert\le 
 2^{-(b-a)^{-1}}2^n (b-a)^{-n}M_n,\]
for any $x\in (a,b)\ss I\bs C$, we can reason as before to conclude that $h\in\cC^\infty(I,\bR)$.

Define $f(x)=\int_0^yh(x)\dif y$. $f\in\cC^\infty(I,\bR)$ and $f$ is strictly increasing. Indeed, since $C$ is totally disconnected, given $x,y\in I$, $x<y$, there exist $t,s\in(x,y)$, $t<s$ such that $[t,s]\ss I\bs C$, so $h(x)>0$ in $[t,s]$ and hence $f(y)-f(x)=\int_x^yh(z)\dif z>0$.

Given a connected component $(a,b)$ of $I\bs C$, \[\mu(f(a,b))=f(b)-f(a)=\int_a^bh(z)\dif z.\]
Thus, given that $f$ is strictly increasing, $I\bs C$ has a countable number of connected components and that the Lebesgue measure is $\sigma$-additive, $\mu(f(I\bs C))=\int_{I\bs C}h(z)\dif z$. $f(I)$ is also measurable because it is an interval. Since $f$ is strictly increasing, $f(C)\cap f(I\bs C)=\emptyset$ and  $f(C)=f(I)\bs f(I\bs C)$ and, therefore,  $f(C)$ is measurable. Hence, \begin{align*}\mu(f(C))=& \mu(f(I))-\mu(f(I\bs C))=f(1)-f(0)-\int_{I\bs C}h(z)\dif z\\ = & \int_0^1h(z)\dif z-\int_{I\bs C}h(z)\dif z=\int_Ch(z)\dif z=0.\end{align*}
Now, since $\mu(C)>0$, there exists  $D\ss C$ such that $D\not\in\cL$ \cite[Exercise~29, pg.~39]{Folland}. We have that $f(D)\ss f(C)$, so $\mu^*(f(D))\le\mu^*(f(C))=0$ and, therefore, $\mu^*(f(D))=0$, so by the completeness of $\mu$, $f(D)\in\cL$. Finally, $f^{-1}(f(D))=D\not\in\cL$, so $f$ is not $\cL/\cL$ measurable.

We now define a  $\cL/\cB$ measurable function $g$ in such a way that $g\circ f$ is not  $\cL/\cB$ measurable. Let $g$ be the characteristic function of the set $f(D)$. $g$ is measurable since $f(D)$ is. Furthermore, $\{1\}\in\cB$, but $(g\circ f)^{-1}(\{1\})=f^{-1}(f(D))=D\not\in\cL$, so $g\circ f$ is not $\cL/\cB$ measurable.
\end{exa}

Example~\ref{exace} has shown that, if we are to provide sufficient conditions for $g\circ f$ to be $\cL/\cL$ measurable, we will have to look further away than the regularity of $f$. In fact, it is clear from Example~\ref{exace} that the behavior of $f$ when it takes inverse images is determinant on the behavior of the composition, so we will try first to impose conditions on $f^{-1}$ in the case $f$ is invertible.

For the next results we consider $\Omega,\Lambda\in\cL$, $f:(\Omega,\cL)\to(\bR,\cB)$, $g:(\Lambda,\cL)\to(\bR,\cB)$, $f(\Omega)\ss\Lambda$.

\begin{lem}\label{lemac} If $g$ is $\cL/\cB$ measurable, $f$ invertible and $f^{-1}$ absolutely continuous, then  $g\circ f$ is $\cL/\cB$ measurable.
\end{lem}
\begin{proof}
    Let $B\in\cB$. Then $g^{-1}(B)\in\cL$. Since $f^{-1}$ is absolutely continuous it takes $\cL$ sets to $\cL$ sets \cite[p. 250]{Natanson}, so $f^{-1}(g^{-1}(B))\in\cL$.
\end{proof}
Observe that Lemma~\ref{lemac} crucially avoids the circumstances of Example~\ref{exace}, since, in that case, the function $f^{-1}$ was not absolutely continuous (since it did not map $\cL$ sets to $\cL$ sets). This illustrates that, in general, even if $f$ is absolutely continuous, $f^{-1}$ needs not to be --see \cite{Spataru2004,Amo2014}. A necessary sufficient condition for $f^{-1}$ to be absolutely continuous (in the case the domain is an interval) can be found in the following result --cf.  \cite[Lemma~2.2]{Cabada2000}, \cite{Spataru2004}.
\begin{lem}\label{critd0}
 Let $J=[c, d]$ be a compact real interval and $f$ a strictly increasing and absolutely continuous function on $J$.
$f^{-1}$ is absolutely continuous if and only if $f^{\prime}(x) \neq 0$ for a.e. $x \in J$.
\end{lem}
Lemma~\ref{critd0} raises the question as to whether the condition $f^{\prime}(x) \neq 0$ a.e. is enough to show that $g\circ f$ is measurable without asking for $f$ to be invertible, which is a very stringent condition. We provide an answer to this question in Corollary~\ref{dna}, but before stating it we will need the following results.
\begin{lem}\label{lemcont}
    Assume $f$ and $g$ are $\cL/\cB$ measurable. If for every $N\ss\bR$ such that $\mu^*(N)= 0$ we have that $f^{-1}(N)\in\cL$, then
    $g\circ f$ is $\cL/\cB$ measurable.
\end{lem}
\begin{proof}
    Let $B\in\cB$. Then $g^{-1}(B)\in\cL$, so, by the regularity of $\mu$, $g^{-1}(B)=C\sqcup N$ where $\mu^*(N)=0$ and $C\in\cB$. By hypothesis,  $f^{-1}(N) \in\cL$. On the other hand, since $C\in\cB$ and $f$ is $\cL/\cB$ measurable, $f^{-1}(C)\in\cL$. Therefore, $f^{-1}(g^{-1}(B))=f^{-1}(C\sqcup N)=f^{-1}(C)\sqcup f^{-1}(N)\in\cL$. We conclude that $g\circ f:(\Omega,\cL)\to(\bR,\cB)$ is measurable.
\end{proof}

Remember that the Lebesgue measure is complete, that is, if $\mu^*(A)=0$ then $A\in\cL$ and $\mu(A)=0$, so we obtain the following corollary.
\begin{cor}\label{corn}
   Assume $f$ and $g$ are $\cL/\cB$ measurable. If for every $N\in\cL$ such that $\mu(N)= 0$ we have that $\mu^*(f^{-1}(N))=0$ then
    $g\circ f$ is $\cL/\cB$ measurable.
\end{cor}
\begin{rem}The hypotheses in Lemma~\ref{lemcont} are sharp in that, if there exists $N\ss\bR$ such that $\mu^*(N)= 0$ and $f^{-1}(N)\not\in\cL$, then there exists a $\cL/\cB$ measurable $g$ such that $g\circ f$ is not $\cL/\cB$ measurable. It is enough to take $g=\chi_{N}$.
\end{rem}
%

\begin{rem}
The condition occurring in Corollary~\ref{corn} is reminiscent of the  Lusin $N$-property:
\begin{equation}\tag{N}\mu(A)= 0\text{ implies }\mu(f(A))=0\text{ for every } A\in\cL,\end{equation}
which is intimately related to absolute continuity, as the following theorem shows.
\end{rem}
\begin{thm}[{\cite[Theorem 7.7]{Saks1964}}]\label{thmac} Let $I$ be a bounded interval and $f:I\to\bR$ continuous. Then $f$ is absolutely continuous if and only if $f$ satisfies satisfies {\rm(N)} on $I$ and $\int_Af'<\infty$ where $A\ss I$ is the set of points $x$ where $f$ is differentiable and $f'(x)\in(0,\infty)$.
	\end{thm}
%

As stated before, the goal now is to move away from the invertibility of $f$. Furthermore, we would like to work in a context where $f$ need not be differentiable, as this conditions is definitely much stronger than mere measurability.  With this aim we try to generalize the concept of the derivative not being zero to non differentiable functions, which leads to the following definition of the sets $D_f(x)$ and $S_f$.
\begin{dfn}
    Let $f:\Omega\ss\bR\to\bR$. For every $x\in \Omega$ we define the set
    \[Df(x):=\bigcap_{\d\in\bR^+}\ol{\conv\left\{\frac{f(y)-f(x)}{y-x}\ :\ y\in \Omega,\ 0<|y-x|<\d\right\}}^{\ol\bR}\]
where $\ol X^{\ol\bR}$ denotes the closure in the extended real line with the compact interval topology and $\conv$ the convex hull. Observe that, if $x\in \Omega \cap \Omega'$, then $Df(x)$ is a nonmepty closed convex subset of $\ol\bR$.   Let 
 \[S_f:=\{x\in \Omega \cap \Omega '\ :\ 0\in Df(x)\}.\]
\end{dfn}
These definitions lead to the main result we want to introduce.
\begin{thm}\label{dna1}If $f$ and $g$ are $\cL/\cB$ measurable and $\mu^*(S_f)=0$, then  $g\circ f$ is $\cL/\cB$ measurable.
\end{thm}
\begin{proof}
    By Corollary~\ref{corn}, it is enough to check that $\mu^*(S_f)=0$ implies that, for every $A\in\cL$, if $\mu^*(A)= 0$ then $\mu^*(f^{-1}(A))=0$. Equivalently, let $\mu^*(f^{-1}(A))\ne0$ and let us prove that $\mu^*(A)\ne 0$. Let $B\ss f^{-1}(A)$ such that $B$ is bounded and $\mu^*( B)>0$. We can assume, without loss of generality, that $B$ contains no isolated points of $\Omega$, since $\Omega\bs \Omega'$ has to be a countable set because the usual topology is second countable, and thus $\mu^*(\Omega\bs \Omega')=0$. This way we guarantee that $D_f(x)\ne \emptyset$ for $x\in B$. Let $C:=f(B)\ss A$.  It will be enough to show that $\mu^*(C)\ne 0$.

 Since $\mu^*(S_f)=0$ and $\mu^*(B)>0$, we have that either $\mu^*(\{x\in B\ :\ Df(x)\ss\bR^+\})>0$ or $\mu^*(\{x\in B\ :\ Df(x)\ss\bR^-\})>0$. Let us assume  the first case and the second will be analogous. Let $E=\{x\in B\ :\ Df(x)\ss\bR^+\}$ and, for $n\in\bN$, 
 \[E_n:=\left\{x\in B\ :\ \frac{f(y)-f(x)}{y-x}>\frac{1}{n}\text{ for all }y\in \Omega,\ 0<|y-x|<\frac{1}{n}\right\}.\]
     Observe that $E=\bigcup_{n\in\bN}E_n$. Indeed, let $x\in E$.  Since $0\not\in Df(x)$ and $Df(x)$ is a closed set, $\inf Df(x)>0$. Let $\e\in(0,\inf Df(x))$. By definition of $Df(x)$, there exists $\d\in(0,\e)$ such that
       \[\frac{f(y)-f(x)}{y-x}>\e\text{ for every }y\in \Omega,\ 0<|y-x|<\d.\]
  Take $n\in\bN$ such that $\d>\frac{1}{n}$. Then,
        \[\frac{f(y)-f(x)}{y-x}>\frac{1}{n}\text{ for every } y\in \Omega,\ 0<|y-x|<\frac{1}{n},\]       
       that is $x\in E_n$.

      On the other hand, if $x\in E_n$ for some $n\in\bN$,
        \[\frac{f(y)-f(x)}{y-x}>\frac{1}{n}\text{ for every } y\in \Omega,\ 0<|y-x|<\frac{1}{n},\]
        and thus, $z\ge\frac{1}{n}$ for every $z\in Df(x)$, so $x\in E$. 
        
    We conclude that  $E=\bigcup_{n\in\bN}E_n$. Hence, there exists $k\in\bN$ such that $\mu^*(E_k)>0$ for, otherwise  --see  \cite[Corollary 12.1.1]{Munroe}, \[\mu^*(E)=\mu^*\(\bigcup_{n\in\bN}E_n\)=\lim_{n\to\infty}\mu^*(E_n)=0,\]
	which would be a contradiction.

    Let $\d\in(0,\mu^*(E_k))$ be fixed arbitrarily. Using a characterization of the exterior measure --see \cite[Lemma~2.2]{Marquez2021}, consider a countable infinite cover of $E_k$ consisting pairwise disjoint intervals $\{[a_n,b_n)\}_\n$  such that \[\mu^*(E_k)+\d>\sum_{\n}(b_n-a_n)\ge \mu^*(E_k)\] and $b_n-a_n<\frac{1}{k}$ for every $n\in\bN$.

Since $E_k\ss B$ is bounded, $\mu^*(E_k)<\infty$, so let $m\in\bN$ be such that \[\sum_{n=m+1}^\infty(b_n-a_n)<\d.\]
Let $X=\bigcup_{n=1}^m[a_n,b_n)$. Then, $\mu^*(E_k\bs X)<\d$ since $\{[a_n,b_n)\}_{n=m+1}^\infty$ is a cover of $E_k\bs X$. Thus,
\[\mu^*(E_k)\le \mu^*(E_k\cap X)+\mu^*(E_k\bs X)< \mu^*(E_k\cap X)+\d.\]
Therefore, $\mu^*(E_k\cap X)>\mu^*(E_k)-\d$. Now, there exists $j\in\{1,\dots, m\}$ such that $\mu^*(E_k\cap [a_j,b_j))>\frac{1}{m}(\mu^*(E_k)-\d)$ for, if we assume otherwise, then, for every $n=1,\dots, m$,
\[\mu^*(E_k\cap [a_n,b_n))<\frac{1}{m}(\mu^*(E_k)-\d),\]
and, therefore,
\[\mu^*(E_k\cap X)\le \sum_{n=1}^m\mu^*(E_k\cap [a_n,b_n))\le\mu^*(E_k)-\d,\]
which is a contradiction. Hence, for every $x,y\in E_k\cap [a_j,b_j)$, we have that $|x-y|<\frac{1}{k}$, so
\[\frac{1}{k}<\frac{f(y)-f(x)}{y-x}.\]
This implies $f$ is strictly increasing on $E_k\cap [a_j,b_j)$. Consider $f^{-1}:f(E_k\cap [a_j,b_j))\to E_k\cap [a_j,b_j)$. We have that, for every $x,y\in f(E_k\cap [a_j,b_j))$,
\[|f^{-1}(y)-f^{-1}(x)|<k|y-x|,\]
so $f^{-1}$ is a Lipschitz function and, therefore,
\begin{align*}0< & \frac{1}{m}(\mu^*(E_k)-\d)<\mu^*(E_k\cap [a_j,b_j))=\mu^*(f^{-1}(f(E_k\cap [a_j,b_j))))\\ \le & k\mu^*(f(E_k\cap [a_j,b_j)))\le k\mu^*(C).\end{align*}
Thus, $\mu^*(C)>0$ and this ends the proof.
\end{proof}

\begin{rem} Returning to Example~\ref{exace}, observe that, $f'=h$ and $h(y)= 0$ for $y\in C$, so $C\ss S_f$ and, since $\mu(C)>0$, we have that $\mu(S_f)>0$, so we cannot apply Theorem~\ref{dna1}.
%
    \end{rem}
\begin{rem}  If $\Omega\ss\bR$ is an interval and $f$ is continuous (in fact, it is enough for $f$ to be a Darboux function) then
	 \[Df(x):=\bigcap_{\d\in\bR^+}\ol{\left\{\frac{f(y)-f(x)}{y-x}\ :\ y\in \Omega,\ 0<|y-x|<\d\right\}}^{\ol\bR},\]
	 which simplifies the calculations.
	\end{rem}

If $f$ is differentiable at $x$, then, by definition of derivative, $Df(x)=\{f'(x)\}$, so we have the following corollary.
\begin{cor}\label{dna}If $f$ and $g$ are $\cL/\cB$ measurable, $f$ is differentiable a.e. and $\mu^*((f')^{-1}(0))=0$, then  $g\circ f:(\Omega,\cL)\to(\bR,\cB)$ is measurable.
\end{cor}

Even if Corollary~\ref{dna} presents a simple condition to check it is clear that it leaves out common and interesting cases. For instance, the theorem cannot be applied in the case where $f$ is a constant function and $\mu(\Omega)>0$. In order to palliate this fact, we will present a corollary of Corollary~\ref{dna} where local constancy is permitted but, in order to do so, we will need to introduce the concept of essentially open set which will allow us to control those places where the function is locally constant.
\begin{dfn} Let $X\ss\bR$. We say that $X$ is \emph{essentially open} if there exists an open set $U\ss\bR$ such that $\mu^*(X\Delta U)=0$ where $X\Delta U:=(U\bs X)\sqcup(X\bs U)$. \end{dfn}
The notion of a set being essentially open can be rephrased in several ways, as the following Lemma shows.
\begin{lem} Let $X\ss\bR$. The following statements are equivalent:
    \begin{enumerate}
        \item $X$ is essentially open.
       \item  $X\in\cL$ and there exist an open set $U\ss\bR$ and $V,W\in\cL$, such that $V\ss U$, $W\ss X\bs U$, $\mu(V)=\mu(W)=0$ and $X=(U\bs V)\sqcup W$.
       \item $X\in\cL$ and there exist a set $Y\ss\bR$ with a countable number of connected components and $V,W\in\cL$ such that $V\ss Y$,  $W\ss X\bs Y$, $\mu(V)=\mu(W)=0$ and $X=(Y\bs V)\sqcup W$.
    \end{enumerate}
\end{lem}
\begin{proof}
(I)$\Ra$(II) Assume $X$ is essentially open and let $U$ be an open set such that  $\mu^*(X\Delta U)=0$.   Since $U\bs X, X\bs U\ss X\Delta U$, we have that $\mu^*(U\bs X)=\mu^*(U\bs X)=0$, so $V:=U\bs X,W:=X\bs U\in \cL$ and $V\cap W=\emptyset$. Finally, since $U\in\cL$, \[X=(X\cap U)\sqcup(X\bs U)=[U\bs(U\bs X)]\sqcup
    (X\bs U)=(U\bs V)\sqcup W\in\cL.\]

(II)$\Ra$(I) If $V,W,U\in\cL$ are such that $U$ is open, $\mu(V)=\mu(W)=0$ and $X=(U\bs V)\sqcup W$, then
\begin{align*}\mu^*(X\Delta U)= & \mu^*((U\bs X)\sqcup(X\bs U))\le\mu^*(U\bs X)+\mu^*(X\bs U)\\ = & \mu^*(U\bs ((U\bs V)\sqcup W))+\mu^*([(U\bs V)\sqcup W]\bs U)\\ = & \mu^*(U\cap V\cap (U\bs W))+\mu^*(W\bs U)\le\mu^*(V)+\mu^*(W)=0.
    \end{align*}
We conclude that $U$ is essentially open.

(II)$\Ra$(III) Since every open subset of $\bR$ has a countable number of connected components, this is evident.

(III)$\Ra$(II) Let $Y\ss\bR$ with a countable number of connected components and $ V,\til W\in\cL$ such that $V\cap \til W=\emptyset$, $\mu(V)=\mu(\til W)=0$ and $X=(Y\bs V)\sqcup\til W$. Let $\{Y_n\}_{n\in\Lambda}$ be the family of connected components of $Y$. The $Y_n$ are intervals. Let $\L_1$ be the set of those indices $k\in\L$ such that $I_k$ is just a point, $\L_2=\L\bs\L_1$. Define $U:=\bigcup_{n\in\Lambda_2}\mathring Y_n$ and \[W:=\til W\cup\(\left[\(\bigcup_{n\in\Lambda_1}Y_n\)\cup\(\bigcup_{n\in\Lambda_2}Y_n\bs \mathring Y_n\)\right]\bs V\).\] $U$ is open and, since $W$ is the union of a zero measure set and a countable set, $W\in\cL$ and $\mu(W)=0$. Finally, $V\cap W=\emptyset$ and
\[X=(Y\bs V)\sqcup\til W=(U\bs V)\sqcup W.\qedhere\]
    \end{proof}


The next result allows us to work with functions $f$ which are not continuous by separating the continuous and discontinuous parts of $f$.
\begin{thm}\label{thmfin} Assume $f$ and $g$ are $\cL/\cB$ measurable and $f$ is a.e. differentiable and can be expressed as $f_1+f_2$ where $f_1$ is absolutely continuous and $f_2(\Omega)$ is countable. If $(f')^{-1}(0)$ is essentially open, then  $g\circ f$ is $\cL/\cB$ measurable.
\end{thm}
\begin{proof}
   Consider an open set $U\ss\bR$ and $V,W\in\cL$ such that $\mu(V)=\mu(W)=0$ and $(f')^{-1}(0)=(U\bs V)\sqcup W$. Let $B\in\cB$. Then $g^{-1}(B)\in\cL$. Now, \begin{align*}f^{-1}(g^{-1}(B))= [f^{-1}(g^{-1}(B))\cap (U\bs V)]\sqcup[f^{-1}(g^{-1}(B))\cap W]\sqcup [f^{-1}(g^{-1}(B))\bs (f')^{-1}(0)].\end{align*}

    Given that $f'\ne 0$ in $\Omega\setminus (f')^{-1}(0)$, in particular $\mu^*(\{x \in \Omega \setminus (f')^{-1}(0):\; f'(x)=0\})=0$, so we have, by Corollary~\ref{dna}, that $f^{-1}(g^{-1}(B))\bs (f')^{-1}(0)\in\cL$.
    Also, $\mu^*(f^{-1}(g^{-1}(B))\cap W)=0$, so $f^{-1}(g^{-1}(B))\cap W \in\cL$. Finally, since $U$ is open, $U=\bigcup_{n\in\Lambda}U_n$ where the $\{U_n\}_{n\in\Lambda}$ is a countable family of pairwise disjoint open intervals.  Assume $x\in (U_n\bs V)\cap f^{-1}(g^{-1}(B))$. Let $y\in U_n$, $x<y$. Since $f'=0$ on $U_n\bs V$ and $\mu(V)=0$, we have that
    \begin{align*}f(y)-f(x)= & f_1(y)-f_1(x)+f_2(y)-f_2(x)=\int_{[x,y]}f_1'(z)\dif z+f_2(y)-f_2(x)\\= & \int_{[x,y]\bs V}f_1'(z)\dif z+f_2(y)-f_2(x)=f_2(y)-f_2(x).\end{align*}
         Therefore, $f$ takes on $U_n$ a countable number of values. $f$ is also $\cL/\cB$ measurable, so
        $f^{-1}(y)\in\cL$ for every $y\in f(\Omega)$. Thus,
        \begin{align*}f^{-1}(g^{-1}(B))\cap (U\bs V)= & \(\bigcup_{y\in g^{-1}(B)}f^{-1}(y)\)\cap \(\bigcup_{n\in\Lambda}U_n\bs V\)= \bigcup_{\substack{y\in g^{-1}(B)\\n\in\Lambda}}f^{-1}(y)\cap (U_n\bs V).\end{align*}
        The number of sets in this last union that is nonempty is countable and, since $U_n,V, f^{-1}(y)\in\cL$, we have that $f^{-1}(g^{-1}(B))\cap (U\bs V)\in\cL$.

        We conclude that $f^{-1}(g^{-1}(B))\in\cL$, so $g\circ f$ is $\cL/\cB$ measurable.
\end{proof}

From now on we will assume that $\Omega$ is an interval. This is no restriction as we can always extend $f$ to an interval by a constant.

Every function of bounded variation $f$ can be decomposed as $f_a+f_s+f_j$ where $f_a$ is absolutely continuous, $f_s$ is continuous and singular (that is $f'_s=0$ a.e.) and $f_j$ is the jump part \cite[Theorem~3.16]{MonSlaTvr18}. The jump part is a countable sum of step functions and contains all of the discontinuities of $f$. It can be expressed as
\[f(x)=\sum_{t<x}\Delta^+f(t)+\sum_{t\le x}\Delta^-f(t)\]
and clearly it has a countable image, so we have the following corollary of Theorem~\ref{thmfin}.
\begin{cor}\label{corfin}Assume $\Omega$ is an interval. If $g$ is measurable, $f$ is of bounded variation, $f_s=0$ and $(f')^{-1}(0)$ is essentially open, then  $g\circ f$ is $\cL/\cB$ measurable.
\end{cor}

\begin{rem} Again, Corollary~\ref{corfin} cannot be applied to Example~\ref{exace}. Since $f'$ is continuous, $(f')^{-1}(0)$ is closed, but in this case it is not essentially open. To see this, first observe that $(f')^{-1}(0)=C\cup Y$ where $Y$ is a countable set formed by the middle points of the connected components of $I\bs C$. Since $\mu(Y)=0$, $(f')^{-1}(0)$ is essentially open if and only if $C$ is. Assume, that $C$ is essentially open and we will arrive to a contradiction. Let  $U\ss\bR$ be an open set  and $V,W\in\cL$, $V\ss U$, such that $V\cap W=\emptyset$, $\mu(V)=\mu(W)=0$ and $C=(U\bs V)\sqcup W$.

    Take $x\in U\bs V$. Then there exist $r\in\bR^+$ such that $(x-r,x+r)\ss U$. Since $x\in C$, $x\in \pd(I\bs C)$, so there exists $y\in (x-r,x+r)\bs C$. Since $I\bs C$ is open, there exists  $s\in\bR^+$ such that $(y-s,y+s)\ss (x-r,x+r)\bs C$. Thus,
    $\mu(C\cap (x-r,x+r))=\mu(((x-r,x+r)\bs V)\sqcup W)=2r$. Therefore, $\mu((x-r,x+r)\bs C)=0$ and, since  $(y-s,y+s)\ss (x-r,x+r)\bs C$, we have that $\mu((y-s,y+s))=0$, which is  a contradiction.

\end{rem}

\section*{Acknowledgments}
The authors were partially supported by Xunta de Galicia, project ED431C 2019/02, and by the Agencia Estatal de Investigaci\'on (AEI) of Spain under grant MTM2016-75140-P, co-financed by the European Community fund ERDF.

 \bibliography{t}
 \bibliographystyle{spmpsciper}
\end{document}